\documentclass[12pt]{article}
\usepackage{amsfonts}
\usepackage{epsfig}
\usepackage{fancyhdr}
\usepackage{color}
\usepackage[latin1]{inputenc}
\usepackage{amsmath}
\usepackage{mathrsfs}
\usepackage{boxedminipage}
\usepackage{thmdefs}
\usepackage{boxedminipage}
\usepackage{natbib}
\usepackage{graphicx}
\usepackage[colorlinks=true,linkcolor=red,citecolor=blue]{hyperref}
\usepackage{amssymb}
\usepackage{dsfont}
\usepackage{fancyhdr}
\usepackage{graphicx}

\numberwithin{equation}{section}

\renewcommand{\cite}{\citet}
\numberwithin{equation}{section}
\begin{document}
\begin{center}
{\Large \bf A new test procedure of independence in Copula models
via $\chi^2$-divergence}
\end{center}
\noindent Salim BOUZEBDA$^{*}$ \& Amor KEZIOU$^{**}$
\\
\noindent$^{*}$L.S.T.A., Universit\'{e} Paris 6. 175 rue du
Chevaleret--8$^{\rm \grave{e}me}$ {\'e}tage, b\^atiment A, Bo\^{i}te 158, 75013 Paris,\\
salim.bouzebda@upmc.fr\\
\noindent$^{**}$Laboratoire de Math\'ematiques (FRE 3111) CNRS,
Universit\'e de Reims Champagne-Ardenne and LSTA-Universit\'e
Paris 6. UFR Sciences, Moulin de la Housse, B.P. 1039,
                   51687 Reims,
                   France.\\
\noindent amor.keziou@upmc.fr \vskip 3mm \noindent Key Words:
dependence function; Gumbel copula; pseudo-likelihood; asymptotic
theory; multivariate rank statistics;
semiparametric estimation. \vskip
3mm

\vskip 3mm
\noindent{\bf Abstract.} We introduce a new test procedure of
independence in the framework of parametric copulas with unknown
marginals. The method is based essentially on the dual
representation of $\chi^2$-divergence on signed finite measures.
The asymptotic properties of the proposed estimate and the test
statistic are studied under the null and  alternative hypotheses, with
simple and standard limit distributions both  when the parameter
is an  interior point or not.

\section{Introduction and motivations}Parametric models for copulas have been intensively investigated
during the last decades. \emph{Copulas} have become of great interest in
applied statistics, because of the fact that they constitute a
flexible and robust way to model dependence between the marginals of
random vectors. The reader may refer to the following books for excellent expositions of the basics of copula theory : \cite{Nelsen1999} and \cite{Joe1997}.
In this framework, semiparametric inference
methods, based on \emph{pseudo-likelihood}, have been applied to
copulas by a number of authors (see, e.g., \cite{Shih_Louis1995},
\cite{Wang_Ding2000}, \cite{tsukahara2005} and the references
therein). Throughout the available literature, investigations on
the asymptotic properties of parametric estimators, as well as the
relevant test statistics, have privileged the case where the
parameter is an interior point of the admissible domain. However,
for most parametric copula models of interest, the boundaries of
the admissible parameter spaces include some important parameter
values, typically among which, that corresponding to the
independence of margins. This paper concentrates on this specific
problem. We aim, namely, to investigate parametric inference
procedures, in the case where the parameter belongs to the
boundary of the admissible domain. In particular,  the usual limit
laws both for parametric copula estimators and test statistics
become invalid under these limiting cases, and, in particular,
under marginal independence. Motivated by this observation, we
will introduce a new semiparametric inference procedure based on
\emph{$\chi^2$-divergence} and \emph{duality} technique. We will
show that the proposed estimator remains asymptotically normal,
even under the marginal independence assumption. This will allow
us to introduce test statistic of independence, his study will
be made, both under the null and alternative hypotheses.\vskip5pt
\noindent It is well known since the work of \cite{Sklar1959} that
the joint behavior of a bivariate vector $(X_{1},X_{2})$ with d.f.
$\mathbf{F}(x_{1},x_{2}):=\mathbf{P}(X_1\leq x_1,X_2\leq x_2)$, and continuous
marginal d.f.'s $F_{i}(x_i) := P(X_{i} \leq x_i),~i=1,2$, is
characterized by the copula (or dependence function) $C(\cdot,\cdot)$
associated with $\mathbf{F}(\cdot,\cdot)$. The copula function is defined,
for all $(u_{1},u_{2})$ $\in (0,1)^2$, through the identity
\begin{gather*}
C(u_{1},u_{2}) := \mathbf{P}\left\{F_{1}(X_{1}) \leq u_{1},F_{2}(X_{2})
\leq u_{2}\right\}.
\end{gather*}
Many useful multivariate models for dependence between $X_1$ and
$X_2$ turn out to be generated by \emph{parametric} families of
copulas of the form $\left\{C_\theta;~\theta\in\Theta\right\}$,
typically indexed by a vector valued parameter
$\theta\in\Theta\subseteq\mathds{R}^p$ (see, e.g.,
\cite{Kimeldorf1075_1}, \cite{Kimeldorf1075_2}, and
\cite{Joe1993}). The nonparametric approach to copula estimation
has been initiated by \cite{deheuvels1979b}, who introduced and
investigated the \emph{empirical copula process}. In addition,
\cite{Deheuvels1981b,deheuvels1981} described the limiting
behavior of this empirical process (see, also
\cite{fermanianradulovicdragan2004} and the references therein).
In this paper, we consider semiparametric copula models with
unknown marginals.\vskip5pt \noindent In order to estimate
the unknown \emph{true} value of the parameter $\theta\in\Theta$,
which we denote, throughout the sequel, by $\theta_T\in\Theta$,
some \emph{semiparametric estimation} procedures, based on the
maximization, on the parameter space $\Theta$, of properly chosen
\emph{pseudo-likelihood} criterion, have been proposed by
\cite{Oakes1994}, and studied by \cite{Genest_Ghoudi_Rivest1995},
\cite{Shih_Louis1995}, \cite{Wang_Ding2000} and
\cite{tsukahara2005} among others. In each of these papers, some
asymptotic normality properties are established for
$\sqrt{n}\big(\tilde{\theta}-\theta_T\big)$, where
$\tilde{\theta}=\tilde{\theta}_n$ denotes a properly chosen
estimator of $\theta_T$. This is achieved, provided that
$\theta_T$ lies in the \emph{interior}, denoted by
$\mathring{\Theta}$, of the parameter space $\Theta\subseteq
\mathds{R}^p$. On the other hand, the case where $\theta_T\in
\partial \Theta:=\overline{\Theta}-\mathring{\Theta}$ is a
\emph{boundary value} of $\Theta$, has not been studied in a
systematical way until present. Moreover, it turns out that, for
the above-mentioned estimators, the asymptotic normality of
$\sqrt{n}\big(\tilde{\theta}-\theta_T\big),$ may fail to hold for
$\theta_T\in
\partial \Theta$; indeed, under some regularity conditions, when $\theta$ is univariate,  we can
prove that the limit law is the distribution of
$Z\mathds{1}_{(Z\geq 0)}$ where $Z$ is a centred normal variable,
and that the limit law of the  generalized pseudo-likelihood ratio
statistic is a mixture of chi-square laws with one degree of
freedom and Dirac measure  at zero; see
\cite{bouzebda-keziou2008}. Furthermore, when the parameter is
multivariate, the derivation of the limit distributions under the
null hypothesis of independence, becomes much more complex; see
\cite{selflieng1987}. Also, the limit distributions are not
standard which yields formidable numerical difficulties to
calculate the critical value of the test. We cite below some
examples of parametric copulas, for which marginal independence is
verified for some specific values of the parameter $\theta$, on
the boundary $\partial \Theta$ of the admissible parameter set
$\Theta$.  We start with examples for which $\theta$ varies within
subsets of $\mathds{R}$. Such is the case for the extreme value copulas, namely
\begin{equation}
C_{A}(u_1,u_2):=\exp \left\{\log u_1u_2A\left(\frac{\log u_1}{\log u_1u_2}\right)\right\},
\end{equation}
where $A(\cdot)$ is a convex function on $[0,1]$,  satisfying
\begin{enumerate}
\item[-]   $A: [0, 1]\mapsto[1/2, 1]$ such that $\max(t, 1-t)\leq
A(t)\leq 1$ for all $0\leq t\leq 1$.
\end{enumerate}
For
\begin{equation}\label{gumbel_copula_f1}
A(t):=A_\theta(t)=(t^\theta+(1-t)^\theta)^{1/\theta} ; ~~\theta \in [1,\infty[
\end{equation}
 we have \cite{Gumbel1960} family of copulas, which is one of the most
popular model used to model bivariate extreme values. For
\begin{equation}\label{3a}
A_\theta(t)=1-(t^{-\theta}+(1-t)^{-\theta})^{-1/\theta}; ~~\theta\in [0, \infty[
\end{equation}
we obtain \cite{Galambos1975} family of copulas. Finally for
\begin{equation}\label{4}
A_\theta(t)=t\Phi\left(\theta^{-1}+\frac{1}{2}\theta\log\left(\frac{t}{1-t}\right)\right)+(1-t)\Phi\left(\theta^{-1}-
\frac{1}{2}\theta\log\left(\frac{t}{1-t}\right)\right),
\end{equation}
where $\theta \in[0,\infty[$ and $\Phi(\cdot)$ denoting the standard normal
$N(0,1)$ distribution function, we obtain the \cite{Husler_Reiss1989} family of copulas.
A useful family of copulas, due to \cite{Joe1993}, is given, for
$0<u_1,u_2<1$, by
\begin{eqnarray}\label{2}
C_\theta(u_1,u_2):=1-\left[(1-u_1)^\theta+(1-u_2)^\theta-(1-u_1)^\theta
(1-u_2)^\theta\right]^{1/\theta}; ~~\theta\in
[1,\infty[.
\end{eqnarray}
The
Gumbel-Barnett copulas are given, for $0<u_1,u_2<1$, by
\begin{gather}\label{5}
C_{\theta}(u_1,u_2) := u_1u_2\exp\left\{-(1-\theta)(\log u_1)(\log
u_2)\right\}; ~~ \theta \in [0,1].
\end{gather}
The Clayton copulas of positive dependence are such that, for $0<u_1,u_2<1$,
\begin{equation}\label{6}
C_\theta(u_1,u_2)=\left(u_1^{-\theta}+u_2^{-\theta}-1\right)^{-1/\theta}; ~~\theta \in ]0,\infty[.
\end{equation}
Parametric families of copulas with parameter $\theta$ varying in
$\mathds{R}^p$, for some $p\geq 2$, include the following
classical examples. Below, we set
$\theta=\big(\theta_1,\theta_2\big)^{\top}\in\mathds{R}^2$.
\begin{eqnarray}\label{7}
&&C_{\theta}(u_1,u_2):=\left\{1+\left[(u_1^{-\theta_1}-1)^
{\theta_2}+(u_2^{-\theta_1}-1)^{\theta_2}\right]^{1/\theta_2}
\right\}^{-1/\theta_1},
~\theta \in ]0,\infty[\times [1,\infty[;\\
\label{8}
&&C_{\theta}(u_1,u_2):=\exp\Big\{-\Big[{\theta_2}^{-1}\log
\Big(\exp\left(-{\theta_2}(\log
u_1)^{\theta_1}\right)\\
&&\qquad+\exp\left(-{\theta_2}(\log u_2)^{\theta_1}\right)
-1\Big)\Big]^{1/\theta_1}\Big\},~\theta \in [1,\infty[\times
]0,\infty[.\nonumber
\end{eqnarray}
For other examples of the kind, we refer to
\cite{Joe1997}.\vskip5pt \noindent For each of the above examples,
the independence case $C_{\theta_T}(u_1,u_2)=u_1u_2$ (or $A(t)=1$) occurs at the
boundary of the parameter space $\Theta$, i.e., when $\theta_T=1$
for the models (\ref{gumbel_copula_f1}), (\ref{2}) and (\ref{5}),
$\theta_T=0$ for the models (\ref{3a}), (\ref{4}) and (\ref{6}),
$\theta_T=(0,1)^\top$ for the bivariate parameter model (\ref{7}), and
$\theta_T=(1,0)^\top$ for the bivariate parameter model (\ref{8}). In
the sequel, we will denote by $\theta_0$ the value of the
parameter (when it exists), corresponding to the independence of
the marginals, i.e., the value of the parameter for which we
have\begin{equation*}C_{\theta_0}(u_1,u_2):=u_1u_2,~ \text{ for
all } (u_1,u_2)\in (0,1)^2.\end{equation*} Hence, $\theta_0=1$ for
the models (\ref{gumbel_copula_f1}), (\ref{2}) and (\ref{5}),
$\theta_0=0$ for the models (\ref{3a}), (\ref{4}) and (\ref{6}),
$\theta_0=(0,1)^\top$ for the model (\ref{7}), and $\theta_0=(1,0)^\top$ for
the model (\ref{8}). Note that for the models (\ref{3a}), (\ref{4}), (\ref{6}), (\ref{7}) and (\ref{8}),
 $C_{\theta_0}(u_1,u_2)=u_1u_2$ is naturally defined to be the limit of $C_{\theta}(\cdot,\cdot)$
 when $\theta$ tends to $\theta_0$ with values in $\Theta$.
 We denote $c_\theta(\cdot,\cdot):=\frac{\partial^2}{\partial
u_1\partial u_2}\,C_\theta(\cdot,\cdot)$ the density of $C_\theta(\cdot,\cdot)$ and
we define $c_{\theta_0}(\cdot,\cdot)$ to be the limit
 of $c_{\theta}(\cdot,\cdot)$
 when $\theta$ tends to $\theta_0$ with values in $\Theta$. Hence, we can show that for all the above models
 $c_{\theta_0}(u_1,u_2)=1$ for all $0<u_1,u_2<1$. \vskip5pt \noindent In contrast with the
preceding examples, where $\theta_0\in\partial\Theta$ is a
boundary value of $\Theta$, the case where $\theta_0$ is an
interior point of $\Theta$ may, at times, occur, but is more
seldom. An example where $\theta_0\in\mathring{\Theta}$ is given
by the Farlie-Gumbel-Morgenstern (FGM) copula, defined by
\begin{equation}\label{FGM}
C_\theta (u_1,u_2):=u_1u_2+\theta u_1
u_2(1-u_1)(1-u_2),~\theta\in\Theta :=[-1,1],
\end{equation}
and for which $\theta_0=0\in\mathring{\Theta}=]-1,1[$.\vskip5pt
\noindent In the present article, we will treat parametric
estimation of $\theta_T$, and tests of the independence assumption
$\theta_T=\theta_0$. We consider both the case where
$\theta_0\in\mathring{\Theta}$ is an interior point of $\Theta$,
and the case where $\theta_0\in \partial \Theta$ is a boundary
value of $\Theta$. To treat this case, we propose a new inference
procedure, based on an estimation of $\chi^2$-divergence by
\emph{duality} technique. This method may be applied independently
of the dimension of the parameter space. Also the limit law of the
estimate of the parameter is normal and the limiting distribution
of the proposed test statistic is $\chi^2$ under independence,
either when $\theta_0$ is an interior point, or when $\theta_0$ is
a boundary point of $\Theta$. The idea is to include the parameter domain
$\Theta$ into an enlarged space, say $\Theta_e$, in order to
render $\theta_0$ an interior point of the new parameter space,
$\Theta_e$. The conclusion is then obtained through an application
of $\chi^2$-divergence and duality technique. Our methods rely on
the fact that, under appropriate assumptions, the definition of
the density $c_\theta(\cdot,\cdot):=\frac{\partial^2}{\partial
u_1\partial u_2}\,C_\theta(\cdot,\cdot)$ of
$C_\theta(\cdot,\cdot)$, pertaining to the models we consider, may
be extended beyond the \emph{standard} domain of variation $\Theta$
of $\theta$. On the other hand, the definition of
$c_\theta(\cdot,\cdot)$ which corresponds to these extensions, is
then, in general, no longer a density, and may, at times, become
negative. For example, such is the case for the parametric models
 (\ref{gumbel_copula_f1}), (\ref{2}), (\ref{3a}) and (\ref{4}),
for which $c_\theta(\cdot,\cdot)$ is meaningful for some
$\theta\not\in\Theta$, but then, becomes negative over some
non-negligible (with respect to Lebesgue's measure) subsets of
$(0,1)^2$. This implies that the log-likelihood of the data is not
properly defined on the whole space $\Theta_e$. For this reason,
we will use the $\chi^2$-divergence between signed finite
measures. We will discuss this problem in more details, below, in
section \ref{fidiv}.\vskip5pt \noindent The remainder of the
present paper is organized as follows. In section \ref{fidiv}, we
present our semiparametric inference procedure, based upon
optimization of the $\chi^2$-divergence between the model
$(C_\theta, \theta\in\Theta_e)$ and the empirical copula
associated to the data, and by using the \emph{dual
representation} of $\chi^2$-divergence. We then derive the
asymptotic limiting distribution of the proposed estimator. It
will become clear later on from our results, that the asymptotic
normality of the estimate holds, even under the independence
assumption, when, either, $\theta_0$ is an interior, or a boundary
point of $\Theta$. The proposed test statistic of independence is also
studied, under the null hypothesis $\mathscr{H}_0$ of
independence, as well as under the alternative hypothesis. The
limiting asymptotic distribution of the test statistic under the
alternative hypothesis is used to derive an approximation to the
power function. An application of the forthcoming results will
allow us to evaluate the sample size necessary to guarantee a
pre-assigned power level, with respect to a specified alternative.
Finally, section \ref{simu} reports a short simulation results, to illustrate the performance of the proposed  test statistic.
The proofs of these results will be postponed to the appendix.
\section{A semiparametric estimation procedure through
$\chi^2$-divergence}\label{fidiv}  As mentioned earlier, the
problem of estimating $\theta$, when $\theta\in\partial\Theta$,
has not been systematically considered in the scientific
literature; and  the classical asymptotic normality property of
the estimators is no longer satisfied. To overcome this
difficulty, in what follows, we enlarge the parameter space
$\Theta$ into a wider space $\Theta_e\supset \Theta$. This is
tailored to let $\theta_0$ become an interior point of $\Theta_e$.
Naturally, we assume that the definition of the function
$c_{\theta}(\cdot,\cdot)$ may be extended to $\Theta_e$. The
difficulty associated with this construction is that, subject to a
proper definition, the densities $c_\theta(\cdot,\cdot)$ of
$C_\theta(\cdot,\cdot)$ with respect to the Lebesgue's measure,
may become negative on some non negligible subsets of $I:=(0,1)^2$
(in this case $c_\theta(\cdot,\cdot)$ becomes the density of a
signed measure, see remark \ref{remark}). Note that just as
Deheuvels's empirical copula is not a copula,
$c_\theta(\cdot,\cdot)$ for $\theta \in \Theta_e$ is not
necessarily a copula density and fail to integrate to $1$. When
such is the case, a semiparametric estimation of $\theta_T$ via
log-likelihood cannot be used. To circumvent this difficulty, we
introduce a new inference procedure, based on $\chi^2$-divergence
method, and duality technique. Recall that the $\chi^2$-divergence
between a bounded signed measure $\mathbf{Q}$, and a probability
$\mathbf{P}$ on $\mathscr{D}$, when $\mathbf{Q}$ is absolutely
continuous with respect to $\mathbf{P}$, is defined by
\begin{equation}\chi^2(Q,P):=\int_{\mathscr{D}}
\varphi\left(\frac{d\mathbf{Q}}{d\mathbf{P}}\right)~d\mathbf{P},\mbox{   where  }
\varphi~:~x\in \mathds{R}\mapsto
\varphi(x):=\frac{1}{2}\left(x-1\right)^2.\end{equation} In the
sequel, we denote by $\chi^2(\theta_0,\theta_T)$ the
$\chi^2$-divergence between $C_{\theta_0}(\cdot,\cdot)$ and $C_{\theta_T}(\cdot,\cdot)$.
Applying the dual representation of $\phi$-divergence obtained by
\cite{BK2005111} Theorem 4.4, we readily obtain that
$\chi^2(\theta_0,\theta_T)$ can be rewritten as
\begin{equation}\label{eqn 1 dual}
\chi^2(\theta_0,\theta_T):=\sup_{f\in\mathscr{F}}\left\{\int_I
f~dC_{\theta_0}-\int_I \varphi^*(f)~dC_{\theta_T}\right\},
\end{equation}
 where
$\varphi^*(\cdot)$ is used to denote the convex conjugate of $\varphi(\cdot)$,
namely, the function defined by
$$\varphi^* ~:~t\in\mathds{R}\mapsto \varphi^*(t):=
\sup_{x\in\mathds{R}}\left\{tx-\varphi(x)\right\}=\frac{t^2}{2}+t,$$
and $\mathscr{F}$ is an arbitrary  class of   measurable
functions, fulfilling the following  conditions \\
$\quad$ $\forall f\in\mathscr{F};~\int |f|~dC_{\theta_0}$ is
finite and
$\varphi'(dC_{\theta_0}/dC_{\theta_T})=\varphi'(1/c_{\theta_T})\in\mathscr{F}$.
Furthermore, the sup in (\ref{eqn 1 dual}) is unique and achieved
at $f=\varphi'(1/c_{\theta_T}).$ \vskip5pt \noindent Define the
new parameter space $\Theta_e$ of $\theta$ as follows. Set
\begin{equation}\label{def de Theta_e }
\Theta_e :=\left\{\theta\in\mathds{R}^d \text{ such that } \int
\left|\varphi'(1/c_\theta(u_1,u_2))\right|~du_1du_2
<\infty\right\}.
\end{equation}
 By choosing the class of functions $\mathscr{F}$, via
$$\mathscr{F}:=\left\{(u_1,u_2)\in I\mapsto
\varphi'(1/c_\theta(u_1,u_2))-1~; ~\theta\in\Theta_e\right\},$$ we
infer from (\ref{eqn 1 dual}) the relation
\begin{eqnarray} \label{eqn 2 dual}
\chi^2(\theta_0,\theta_T):=\sup_{\theta\in\Theta_e}\Big\{\int_{I}
\left(\frac{1}{c_\theta(u_1,u_2)}-1\right)du_1du_2-\int_{I}
\left(\frac{1}{2}\frac{1}{c_\theta(u_1,u_2)^2}-\frac{1}{2}\right)dC_{\theta_T}(u_1,u_2)\Big\}.
\end{eqnarray}
It turns out that the supremum in $(\ref{eqn 2 dual})$  is reached
iff $\theta=\theta_T$. Moreover, in general, $\theta_T$ is an
interior point of the new parameter space $\Theta_e$, especially
under the null hypothesis of independence, namely, when
$\theta_T=\theta_0.$ Set
\begin{equation}
m(\theta,
u_1,u_2):=\int_{I}\left(\frac{1}{c_{\theta}(u_1,u_2)}-1\right)~du_1du_2
 -\left\{
 \frac{1}{2}\frac{1}{{c_{\theta}(u_1,u_2)}^2}-\frac{1}{2}\right\}.
\end{equation}
Consider a random sample $\{(X_{1k},X_{2k}); k = 1,\ldots, n\}$ from the
distribution of $(X_1,X_2)$ denoted by $\mathbf{F}_{\theta_T} (x_1, x_2) := C_{\theta_T} (F_1(x_1), F_2(x_2)) = \mathbf{P}(X_1 \leq x_1,X_2 \leq  x_2).$
In what follows, we propose to estimate the $\chi^2$-divergence
$\chi^2(\theta_0,\theta_T)$ between $C_{\theta_0}(\cdot,\cdot)$ and
$C_{\theta_T}(\cdot,\cdot)$, by
\begin{equation}\label{estim chi2}
\widehat{\chi}^2(\theta_0,\theta_T) := \sup_{\theta \in \Theta_e}
\int_{I} m(\theta,u_1,u_2)~dC_{n}(u_1,u_2),
\end{equation}
and to estimate the parameter $\theta_T$ by
\begin{equation}\label{estimate}
\widehat{\theta}_n := \arg \sup_{\theta \in \Theta_e}
\left\{\int_{I} m(\theta,u_1,u_2)~dC_{n}(u_1,u_2)\right\},
\end{equation}
where  $C_n(\cdot,\cdot)$ is the modified empirical copula, defined by
\begin{equation}\label{ec}
C_n(u_1,u_2):=\frac{1}{n}\sum_{k=1}^{n}\mathds{1}_{\left\{F_{1n}(X_{1k})\leq
u_1\right\}} \mathds{1}_{\left\{F_{2n}(X_{2k})\leq
u_2\right\}},~(u_{1},u_{2})\in I,
\end{equation}
where
$$F_{jn}(t):= \frac{1}{n}\sum_{k=1}^n
\mathds{1}_{]-\infty,t]}(X_{jk}),~j=1,2,$$
and  $\mathds{1}_A$ stands for the indicator function of the event $A$.
\begin{remark}
The choice of the $\chi^2$-divergence among various divergences is motivated by the following statements:
\begin{enumerate}
\item[-] Recall that the $\phi$-divergences between a bounded
signed measure $\mathbf{Q}$, and a probability $\mathbf{P}$ on $\mathscr{D}$, when $\mathbf{Q}$ is absolutely
continuous with respect to $\mathbf{P}$, is defined
by
$$D_\phi(\mathbf{Q},\mathbf{P}):=\int_{\mathscr{D}} \phi\left(\frac{d\mathbf{Q}}{d\mathbf{P}}\right)~d\mathbf{P},$$
where $\phi$ is a proper closed convex function from
$]-\infty,\infty[$ to $[0,\infty[$ with $\phi(1)=0$ and such that
the domain ${\rm dom}\phi:=\{x \in \mathds{R}: \phi(x)<\infty\}$
is an interval with end points $a_{\phi}<1<b_{\phi}$.The
Kullback-Leibler, modified Kullback-Leibler, $\chi^2$, modified
$\chi^2$, and Hellinger  divergences are examples of
$\phi$-divergences; they are obtained respectively for
$\phi(x)=x\log x-x+1$, $\phi(x)=-\log x+x-1$,
$\phi(x)=\frac{1}{2}(x-1)^2$,
$\phi(x)=\frac{1}{2}\frac{(x-1)^2}{x}$, and $\phi(x)=2(\sqrt{x}-1)^2$. We observe  that
Kullback-Leibler, modified Kullback-Leibler, modified
$\chi^2$, and Hellinger  divergences  are infinite when $d\mathbf{Q}/d\mathbf{P}$ takes negative values on non negligible  (with respect to $\mathbf{P}$) subset of $\mathscr{D}$, since the corresponding
$\phi(\cdot)$ is infinite on $(-\infty,0)$. This problem does not hold in the case of $\chi^2$-divergence, indeed the corresponding $\phi(\cdot)$ is finite on $\mathds{R}$.
\item[-] We give an example of copulas for which the likelihood-based procedure fails. We consider the Gumbel copulas $C_\theta(\cdot,\cdot)$ given in (\ref{gumbel_copula_f1}), it's corresponding density copula is defined by
    \begin{eqnarray}
    c_\theta(u_1,u_2):=C_\theta(u_1,u_2)(u_1u_2)^{-1}\frac{(\widetilde{u_1}\widetilde{u_2})^{(\theta-1)}}
    {(\widetilde{u_1}^\theta+\widetilde{u_2}^\theta)^{(2-1/\theta)}}
    \left[(\widetilde{u_1}^\theta+\widetilde{u_2}^\theta)^{(1/\theta)}+\theta-1\right],
    \end{eqnarray}
    where $\widetilde{x}=-\log x$.
We can show that $c_\theta(\cdot,\cdot)$ may takes negative values for some $\theta\in \Theta_e$. In fact $c_{0.7}(u_1,u_2)$ is negative for $(u_1,u_2)\in[0.9,1]^2$, hence the likelihood function is not well defined.
\end{enumerate}
\end{remark}
\begin{remark} The set $\Theta_e$ defined in (\ref{def de Theta_e }) is generally with
non empty interior $\mathring{\Theta}_e$. In particular, we may
check that $\theta_0$ (the value corresponding to independence)
belongs to $\mathring{\Theta}_e$, since the integral in (\ref{def
de Theta_e }) is  finite;  it is equal to zero when
$\theta=\theta_0$, for any  copula density $c_\theta(\cdot,\cdot)$. For
example, in the case of FGM copulas (\ref{FGM}), it is easy to
show that $\Theta_e=\mathds{R}$. However,
it is hard to determine the whole set $\Theta_e$ for some copulas, but in order to test the independence, we need only to prove the existence of a
neighborhood $N(\theta_0)$ of $\theta_0$ for which the integral in (\ref{def de Theta_e }) is finite since we calculate the estimate
 $\widehat{\theta}_n$ in (\ref{estimate}) by Newton-Raphson algorithm using $\theta_0$ as initial point. The explicit calculation
of the integral in (\ref{def de Theta_e }) may be complicated for
some copulas, in such cases we use the Monte Carlo method to
compute this integral.
\end{remark}
\noindent Statistics of the form
\[\Psi_n :=\int_I\psi(u_1,u_2)~dC_{n}(u_1 ,u_2),\]  belong to the general
class of \emph{multivariate rank statistics}. Their asymptotic
properties have been investigated at length by a number of
authors, among whom we may cite \cite{Ruymgaart_Shorack_Zwet1972},
\cite{Ruymgaart1974} and \cite{Ruschendorf1976}. In particular,
the previous authors have provided regularity conditions, imposed
on $\psi(\cdot,\cdot)$, which imply the asymptotic normality of
$\Psi_n$. The corresponding arguments have been modified by
\cite{Genest_Ghoudi_Rivest1995}, \cite{tsukahara2005} and
\cite{fermanianradulovicdragan2004report}, as to establish almost
sure convergence of the estimators they consider. In the same
spirit, the limiting behavior, as $n$ tends to the infinity, of our
estimator  and  test statistic $\mathbf{T}_n$, will make an
instrumental use of the general theory of multivariate rank
statistics. Using some similar arguments as in
\cite{qinlawless1994}, the existence and consistency of our
estimator and
 test statistic will be established through an application of the law of the
iterated logarithm for empirical copula processes, in combination
with general arguments from multivariate rank statistics theory.
In the sequel, without loss of generality, we will limit ourselves
to the case where the parameter is univariate. The extensions of
our result in a multivariate framework may be achieved, at the
price of additional technicalities, and under similar assumptions.
We will make use of the following definitions.
   \begin{definition}
\begin{enumerate}
\item[(i)] Let $\mathscr{Q}$ be the set of continuous functions
$q$ on $[0,1]$ which are positive on $(0,1)$, symmetric about
$1/2$, increasing on $[0,1/2]$ and satisfy
$\int_0^1\{q(t)\}^{-2}dt< \infty$.
 \item[(ii)] A function
$r: (0,1)\longrightarrow (0,\infty)$  is called u-shaped if it is
symmetric about $1/2$ and increasing on $(0,1/2]$. \item[(iii)]
For $0<\beta < 1$ and u-shaped function $r$, we define
$$
r_\beta(t):=\left\{
\begin{array}{ccc}
  r(\beta t)& if & 0<t\leq 1/2; \\
   r\{1-\beta(1-t)\}& if & 1/2<t\leq 1. \\
\end{array}\right.
$$
If for $\beta > 0$ in a neighborhood of $0$, there exists a
constant $M_\beta$, such that $r_\beta \leq M_\beta r$ on $(0,1)$,
then $r$ is called a reproducing u-shaped function. We denote by
$\mathscr{R}$ the set of reproducing u-shaped functions.
\end{enumerate}
\end{definition}\noindent Typical examples of elements in $\mathscr{Q}$ and
$\mathscr{R}$ are given by \begin{equation*}
q(t)=\left[t(1-t)\right]^\zeta,~
0<\zeta<1/2,~~~~r(t)=\varrho\left[t(1-t)\right]^{-\varsigma},~\varsigma\geq
0,~\varrho\geq 0.
\end{equation*}
We will describe the asymptotic properties of the proposed
estimate $\widehat{\theta}_n$ under the following conditions.
\begin{enumerate}
\item [(C.1)] There exist functions $r_{1,k}, r_{2,k} \in
\mathscr{R}$ such that
\begin{equation*}
\Big|\frac{\partial^k}{\partial\theta^k}\,m(\theta_T,u_1,u_2)\Big|\leq
r_{1,k}(u_1)r_{2,k}(u_2),~~\mbox{for}~~k=1,2,
\end{equation*}and\begin{equation*}
\Big|\frac{\partial^3}{\partial\theta^3}\,m(\theta,u_1,u_2)\Big|\leq
r_{1,3}(u_1)r_{2,3}(u_2) \text{ on a neighborhood }  N(\theta_T)
\text{ of } \theta_T,
\end{equation*}
where $\int_{I}\left\{r_{1,k}(u_1)r_{2,k}(u_2)\right\}^2
dC_{\theta_T}(u_1,u_2)< \infty$ for $k=1,2,3;$
    \item[(C.2)]
 The function $(u_1,u_2)\in I\mapsto
\frac{\partial}{\partial\theta}\,m(\theta_T,u_1,u_2)$ is of
bounded variation on $I$;
        \item[(C.3)] For each $\theta$, the function $\frac{\partial}{\partial\theta}\,m(\theta,u_1,u_2): I\longrightarrow
    \mathds{R}$ is continuously differentiable, and there exist functions $r_i \in
    \mathscr{R}$, $\widetilde{r}_i \in
    \mathscr{R}$
    and $q_i \in \mathscr{Q}$ for  $i=1,2$, such that
\begin{equation*}
    \Big|\frac{\partial^2}{\partial\theta\partial u_i}\,m(\theta,u_1,u_2)\Big|\leq
    \widetilde{r}_i(u_i)r_j(u_j), i,j=1,2~ \mbox{and}~ i\neq j,
    \end{equation*}and
\begin{equation*}
\int_{I}\left\{q_i(u_i)\widetilde{r}_i(u_i)r_j(u_j)\right\}
dC_{\theta_T}(u_1,u_2)< \infty,~ i,j=1,2 ~ \mbox{and}~ i\neq j.
\end{equation*}
\end{enumerate}
\noindent Set
\begin{equation}\label{beta1}
\Xi := -
E\left[{\textstyle\frac{\partial^2}{\partial\theta^2}}m\left(\theta_T,F_{1}
 (X_{1k}),F_{2}(X_{2k})\right)\right],
\end{equation}
and
 \begin{equation}\label{sigma33}
 \Sigma^{2} := {\rm Var}\left[{\textstyle\frac{\partial}{\partial\theta}}m\left(\theta_T,F_{1}
 (X_{1}),F_{2}(X_{2})\right) + W_{1}(\theta_T,X_1) + W_{2}(\theta_T,X_2)\right],
\end{equation}
where
\begin{equation*}
 W_{i}(\theta_T,X_i) := \int_{I} \mathds{1}_{\{F_{i}(X_i) \leq
 u_i\}}{\textstyle\frac{\partial^2}{\partial\theta \partial
 u_i}}m\left(\theta_T,u_1,u_2\right)
 c_{\theta_T}(u_1,u_2)~du_1du_2,~i=1,2.
\end{equation*}
We can see that $\Xi$ and $W_{i}(X_i)$ can be defined,
respectively, by
$$\Xi =  E\left[\left({\textstyle\frac{\partial}{\partial\theta}}m\left(\theta_T,F_{1}
 (X_{1k}),F_{2}(X_{2k})\right)\right)^2\right]
$$
and, for $i=1,2$,
$$W_{i}(\theta_T,X_i) =- \int_{I} \mathds{1}_{\{F_{i}(X_i) \leq
 u_i\}}{\textstyle\frac{\partial}{\partial\theta}}m\left(\theta_T,u_1,u_2\right)
{\textstyle\frac{\partial}{\partial
u_i}}m\left(\theta_T,u_1,u_2\right)
 c_{\theta_T}(u_1,u_2)~du_1du_2.
$$
The following theorem describes the asymptotic behavior of the
estimate $\widehat{\theta}_n$ given in (\ref{estimate}).
From now on, $\stackrel{d}{\longrightarrow}$ denotes the convergence
in distribution.
\begin{theorem}\label{th}
Assume that the conditions C.1-C.3 hold.
\begin{enumerate}
\item [(a)] Let $B(\theta_{T},n^{-1/3}):=\left\{\theta \in
\Theta_e, |\theta -\theta_{T} | \leq n^{-1/3} \right\}$. Then, as
$n\rightarrow\infty$, with probability one, the function
$\theta\mapsto\int m\left(\theta,u_1,u_2\right)~dC_{n}(u_1,u_2)$
reaches its
 maximum value at some point $\widehat{\theta}_n$ in the interior of the
 interval $B(\theta_{T},n^{-1/3})$. As a consequence, the estimate
 is consistent almost surely and  satisfies
  \begin{equation*}
  \int_{I} {\textstyle\frac{\partial}{\partial\theta}}
  m\left(\widehat{\theta}_n,u_1,u_2\right)~dC_{n}(u_1 ,u_2)
= 0.
 \end{equation*}
\item [(b)]  As $n\rightarrow\infty$,
\[\sqrt{n}(\widehat{\theta}_n -\theta_T)\stackrel{d}{\longrightarrow}
N(0,\Sigma^2/\Xi^2).\]
\end{enumerate}
\end{theorem}\noindent The proof of Theorem \ref{th} is postponed to section \ref{thhh1}.\section{A test based on ``$\chi^2$-divergence''} One of the
motivations of the present work is to build a statistical test of
independence, based on $\chi^2$-divergence. In the framework of
the parametric copula model, the null hypothesis, namely, the
independence case $C_{\theta}(u_1,u_2) = u_1u_2$ corresponds to
the condition that $\mathscr{H}_0~:~\theta_T = \theta_0$.  We
consider the alternative composite hypothesis
$\mathscr{H}_1~:~\theta_T \neq \theta_0$. The corresponding
generalized pseudo-likelihood ratio statistic is then given by
\begin{gather*}
{\bf
S}_n(\theta_0,\tilde{\theta})=2\log\frac{\sup_{\theta\in\Theta}
\prod_{k=1}^{n}c_{\theta}(\widehat{F}_{1n}(X_{1k}),\widehat{F}_{2n}(X_{2k}))}
{\prod_{k=1}^{n}c_{\theta_0}(\widehat{F}_{1n}(X_{1k}),\widehat{F}_{2n}(X_{2k}))},
\end{gather*}
where, for $j=1,2$,  $\widehat{F}_{jn}$ stands for $n/(n+1)$ times
the marginal empirical distribution function of the $j$-th
variable $X_j$. The rescaling by the factor $n/(n+1)$, avoids
difficulties arising from potential unboundedness of $\log
c_{\theta}(u_1,u_2)$, when either $u_1$ or $u_2$ tends to one.
Since, $\theta_0$ is a boundary value of the parameter space
$\Theta$, we can see, that the convergence in distribution of
${\bf S}_n$ to a $\chi^2$ random variable is not likely to hold.
In order to bring a solution to this problem, we introduce the
statistic
\begin{gather}\label{mmm}
{\bf T}_n(\theta_0,\widehat{\theta}_n):=2n\sup_{\theta \in
\Theta_e} \int_{I} m(\theta,u_1,u_2)~dC_{n}(u_1,u_2).
\end{gather}
Below, we will show that, under the null hypothesis
$\mathscr{H}_0$, the just-given statistic ${\bf
T}_n(\theta_0,\widehat{\theta}_n)$ converges in distribution to a
$\chi^2$ random variable. This property allows us to build a test
of $\mathscr{H}_0$ against $\mathscr{H}_1$, asymptotically of
level $\alpha$. The limit law of ${\bf
T}_n(\theta_0,\widehat{\theta}_n)$ will also be given under the
alternative hypothesis $\mathscr{H}_1$. The following additional
conditions will be needed for the statement of our results.
\begin{enumerate}
\item [(C.4)] We have $$\lim_{\theta \rightarrow
\theta_0}\frac{\partial}{\partial u_i}\,m(\theta,u_1,u_2)=0,$$ and
there exist $M_1>0$ and $\delta_1>0$ such that, for all $\theta$
in some neighborhood of $\theta_0$, one has, for $i=1,2$,
\begin{equation*}
\left|\frac{\partial}{\partial\theta}\,m(\theta,u_1,u_2)\frac{\partial}{\partial
u_i}\,m(\theta,u_1,u_2)c_\theta(u_1,u_2)\right|<
M_1r(u_i)^{-1.5+\delta_1}r(u_{3-i})^{0.5+\delta_1},
\end{equation*}
where $r(u) := u(1-u)$ for $u\in
    (0,1)$;
 \item [(C.5)] There exist a neighborhood $N(\theta_T)$ of $\theta_T$, and  functions $r_i \in \mathscr{R}$   such that
for all $\theta \in N(\theta_T)$, we have
    \begin{equation*}
    \Big|\frac{\partial}{\partial \theta}m(\theta,u_1,u_2)\Big|\leq r_1(u_1)r_2(u_2)~~
\mbox{with}~~\int_{I}\left\{r_{1}(u_1)r_{2}(u_2)\right\}^2
dC_{\theta_T}(u_1,u_2)< \infty;\end{equation*}
 \item [(C.6)] There
exist functions $r_i, \widetilde{r}_i \in \mathscr{R}$, $q_i \in
\mathscr{Q}$, $i=1,2$ such that
\begin{equation*}
\Big|\,m(\theta_T,u_1,u_2)\Big|\leq r_1(u_1)r_{2}(u_2),
\end{equation*}
\begin{equation*}
    \Big|\frac{\partial}{\partial u_i}\,m(\theta,u_1,u_2)\Big|\leq
    \widetilde{r}_i(u_i)r_j(u_j), i,j=1,2~ \mbox{and}~ i\neq j,
    \end{equation*} with
$$\int_{I}\left\{r_{1}(u_1)r_{2}(u_2)\right\}^2
dC_{\theta_T}(u_1,u_2)< \infty,$$ and\begin{equation*}
\int_{I}\left\{q_i(u_i)\widetilde{r}_i(u_i)r_j(u_j)\right\}
dC_{\theta_T}(u_1,u_2)< \infty,~ i,j=1,2 ~ \mbox{and}~ i\neq j.
\end{equation*}\end{enumerate}
\begin{theorem}\label{thh}
Assume that conditions C.1-C.4 hold. Then, under the null
hypothesis $\mathscr{H}_0$, the statistic ${\bf T}_n$ converges in
distribution to a $\chi_{1}^2$ random variable (with $1$ degree of
freedom).
\end{theorem}
\noindent The proof of Theorem \ref{thh} is postponed until section \ref{thhh2}.
\begin{theorem}\label{thhll}
 Assume that conditions C.1-C.3 and C.5-C.6 hold. Then,
under the alternative hypothesis  $\mathscr{H}_1$, we have
\begin{equation*}
\sqrt{n}\left(\frac{{\bf
T}_n}{2n}-\chi^2(\theta_0,\theta_T)\right)
\end{equation*}
converges to a centered normal random variable with variance
\begin{equation}\label{k1}
\sigma^2_{\chi^2} := {\rm Var}\left[m\left(\theta_T,F_{1}
 (X_{1k}),F_{2}(X_{2k})\right) + Y_{1}(\theta_T,X_1) + Y_{2}(\theta_T,X_2)\right],
\end{equation}
where
\begin{equation*}
 Y_{i}(\theta_T,X_i) := \int_{I} \mathds{1}_{\{F_{i}(X_i) \leq
 u_i\}} {\textstyle\frac{\partial}{\partial u_i}}m\left(\theta_T,u_1,u_2\right)c_{\theta_T}(u_1
,u_2)~du_1du_2.
\end{equation*}
\end{theorem}
\noindent The proof of Theorem \ref{thhll} is postponed until section \ref{thhh2}.\\
\begin{remark} An application of Theorem \ref{thh}, leads to reject the
null hypothesis $\mathscr{H}_0:\theta_T=\theta_0$, whenever the
value of the statistic ${\bf T}_n$ exceeds $q_{1-\alpha}$, namely,
the $(1-\alpha)$-quantile of the $\chi^2_1$ law. The test corresponding to this rejection rule is then,
asymptotically of level $\alpha$, when $n\rightarrow\infty$.
Accordingly, the critical region is given by
\begin{equation}
CR := \left\{{\bf T}_n > q_{1-\alpha}\right\}.
\end{equation}
The fact that this test is consistent follows from Theorem \ref{thhll}.
Further, this theorem can be used to give an approximation to the
power function $\theta_T\mapsto
\beta(\theta_T):=P_{\theta_T}\left\{CR\right\}$ in a similar way
to \cite{KeziouLeoni2005} and \cite{KeziouLeoni2007}. We so obtain
that
\begin{equation}\label{marre2}
\beta(\theta_T)\approx 1-\Phi\left(
\frac{\sqrt{n}}{\sigma_{\chi^2}}
\left(\frac{q_{1-\alpha}}{2n}-\chi^2(\theta_0,\theta_T)\right)\right),
\end{equation}
where $\Phi$ denotes, as usual, the distribution function of
$N(0,1)$ standard normal random variable. A useful consequence of
(\ref{marre2}) is the possibility of computing an approximate
value of the sample size ensuring a specified power
$\beta(\theta_T)$, with respect to some pre-assigned alternative
$\theta_T\neq\theta_0.$ Let $n_0$ be the positive root of the
equation
\begin{equation}
\beta = 1-\Phi\left( \frac{\sqrt{n}}{\sigma_{\chi^2}}
\left(\frac{q_{1-\alpha}}{2n}-\chi^2(\theta_0,\theta_T)\right)\right),
\end{equation}
which can be rewritten as
\begin{equation}
n_0=\frac{(a+b)- \sqrt{a(a+2b)}}{2\chi^2(\theta_0,\theta_T)},
\end{equation}
where $a := \sigma^2_{\chi^2}\left(\Phi^{-1}(1-\beta)\right)^2$
and $b:=q_{1-\alpha} \chi^2(\theta_0,\theta_T)$. The sought-after
approximate value of the sample size is then given by
$$n^*:=\lfloor n_0\rfloor+1,$$
where $\lfloor u\rfloor$ denote the integer part of $u$.
\end{remark}\begin{remark}\label{remark}
If the parameter $\theta$ does not belong to $\Theta$, i.e.,
$\theta \in \Theta_e-\Theta$, the densities $c_\theta(\cdot,\cdot)$, with
respect to the Lebesgue's measure, associated to $C_\theta(\cdot,\cdot)$, may
become negative on some non negligible subsets of $I$. So, in this
case
\begin{equation}\label{ff}
C_\theta(u_1,u_2):=\int_{(0,u_1)\times(0,u_2)}c_\theta(u_1,u_2)
~d\lambda(u_1,u_2);~\forall u\in I,
\end{equation}
may be a signed measure (and not a copula) for some
$\theta\in\Theta_e-\Theta$. Hence, the formula (\ref{estimate})
may lead to a value of the estimator $\widehat{\theta}_n$
belonging to $\Theta_e-\Theta$, and the corresponding
$C_{\widehat{\theta}_n}$ is not necessarily a copula. So, for
point estimation, the pseudo-maximum likelihood estimator
(restricted to vary in the admissible domain) should be used
instead of (\ref{estimate}). The latter estimator may not be
meaningful, and is likely to have a larger mean square error. The
main advantage of our formula (\ref{estim chi2}), where the
supremum is taken on the extended space $\Theta_e$, instead of the
admissible domain $\Theta$, is that it permits easily  to build a
test of independence even when the dimension of the parameter
space $\Theta$ is larger than one. We can show that the proposed
test statistic, based on the formula (\ref{estim chi2}), has a
$\chi^2(p)$ limit law with $p=dim(\Theta)$ degrees of freedom.
\end{remark}
\begin{remark}
The asymptotic variances (\ref{sigma33}) and (\ref{k1})  may be
consistently estimated respectively by the sample variances of
\begin{equation}
 {\textstyle\frac{\partial}{\partial\theta}}\,m\left(\widehat{\theta}_n,F_{1n}
 (X_{1,k}),F_{2n}(X_{2,k})\right) + W_{1}(\widehat{\theta}_n,X_{1,k}) +
 W_{2}(\widehat{\theta}_n,X_{2,k}), \quad k=1,\ldots,n,
\end{equation}
\begin{equation}
 m\left(\widehat{\theta}_n,F_{1n}
 (X_{1,k}),F_{2n}(X_{2,k})\right) + Y_{1}(\widehat{\theta}_n,X_{1,k}) +
 Y_{2}(\widehat{\theta}_n,X_{2,k}), \quad k=1,\ldots,n,
\end{equation}
as was done in \cite{Genest_Ghoudi_Rivest1995}. Similarly, the
parameter $\Xi$ in (\ref{beta1}) may be consistently estimated by
the sample mean of
\begin{equation}
 \left[{\textstyle\frac{\partial}{\partial\theta}}\,m\left(\widehat{\theta}_n,F_{1n}
 (X_{1,k}),F_{2n}(X_{2,k})\right)\right]^2, \quad
 k=1,\ldots,n.
\end{equation}
\end{remark}
\section{Simulation results}\label{simu}
In this section, we present some simulation results aiming  to
illustrate the theoretical results of Theorem \ref{thh}
and Theorem \ref{thhll}.\subsection{The chi square approximation}
We illustrate the accuracy of the approximation of the statistic $\mathbf{T}_n$ by its limit law $\chi_1^2$  under the null hypothesis of independence of
marginals; see Theorem
\ref{thh}.  We consider the Clayton, FGM and Gumbel copulas
 given in (\ref{6}), (\ref{FGM}) and (\ref{gumbel_copula_f1}) respectively.
 We use the Q-Q plot of the empirical quantiles of the proposed statistic
 $\mathbf{T}_n$  versus the quantiles of the $\chi_1^2$ law. In Figure 1, the  Q-Q plots are obtained from 1000
 independent runs of samples with sizes $n=100$ and  $n=500$ for Clayton and FGM copulas, and from 500
 independent runs of samples with sizes $n=100$ and  $n=500$ for Gumbel copula. We
observe that the approximation
is  good even for moderate sample sizes. The integral in the expression of
$m(\theta,\cdot,\cdot)$ is calculated by the Monte Carlo method,
and the supremum in (\ref{mmm}) is considered on the extended
space $\Theta_e$;  it has been computed on
a neighborhood of $\theta_0:=0$ by the Newton-Raphson algorithm
taking $\theta_0:=0$ as an initial point.

\begin{figure}[!h]\label{figd1}
\begin{center}
\begin{boxedminipage}[t]{13cm}
\begin{center}
\begin{tabular}{cc}
\mbox{\epsfig{file=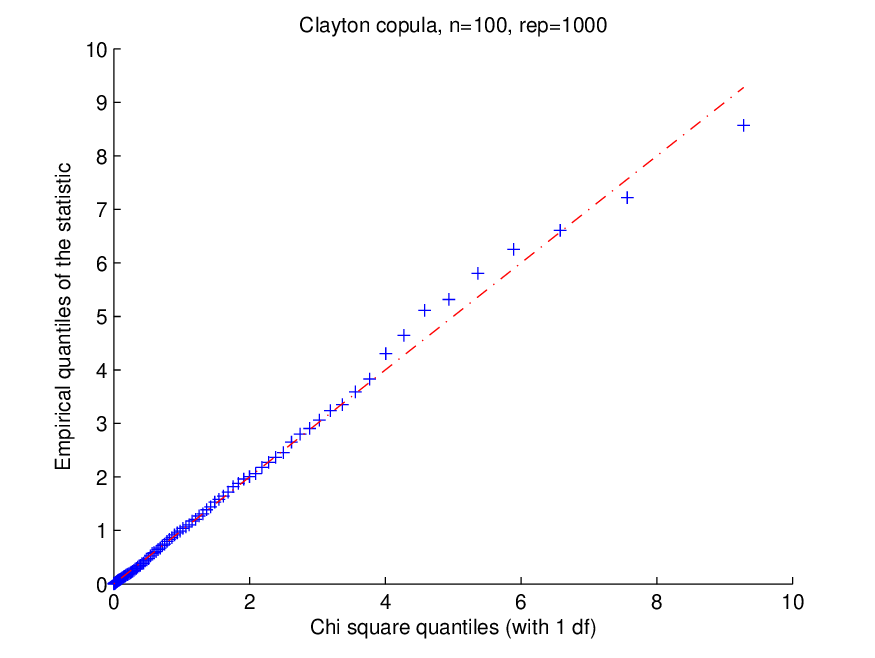,width=5cm,height=4cm}} &
\mbox{\epsfig{file=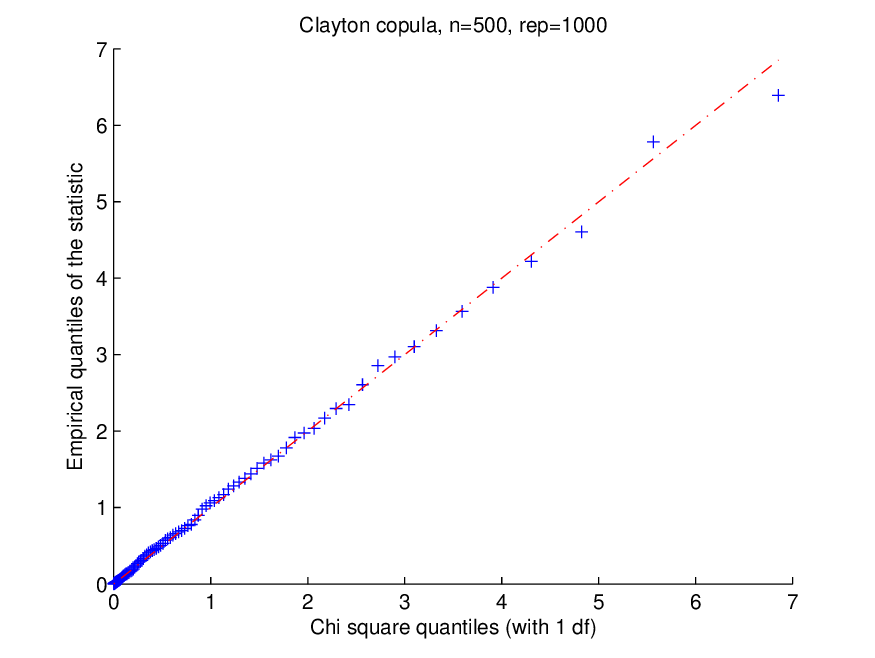,width=5cm,height=4cm}}\\
\mbox{\epsfig{file=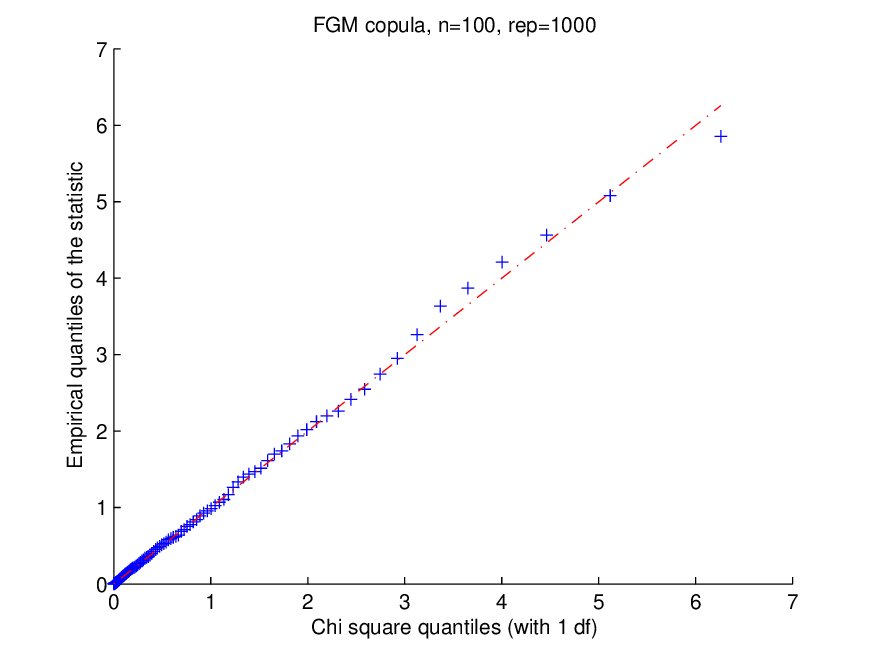,width=5cm,height=4cm}}&
\mbox{\epsfig{file=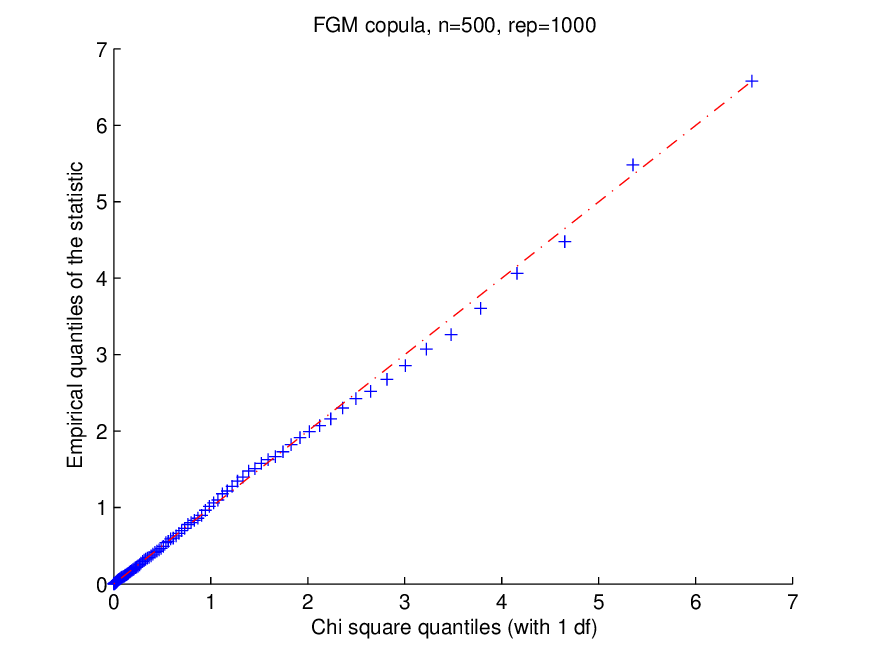,width=5cm,height=4cm}}\\
\mbox{\epsfig{file=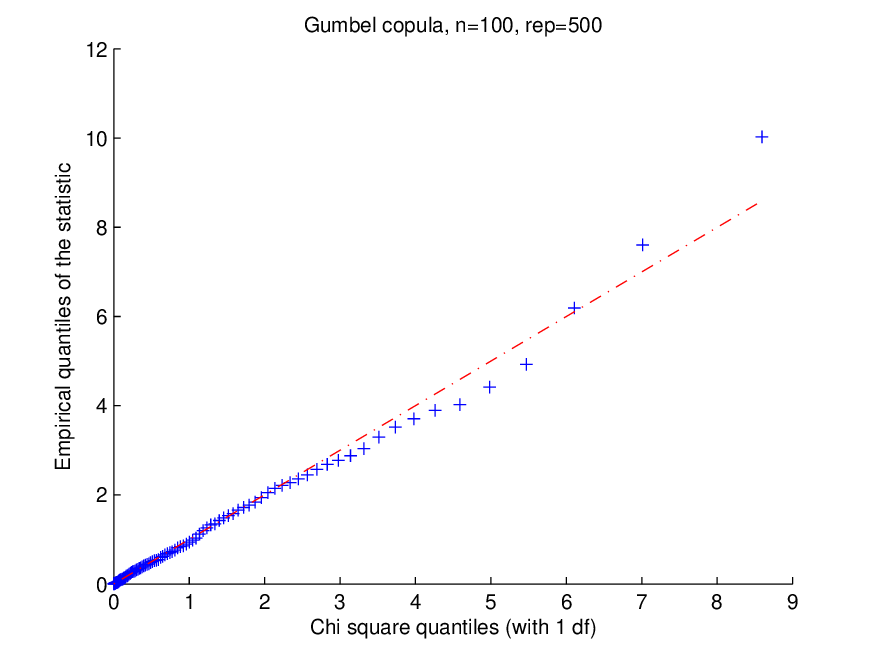,width=5cm,height=4cm}}&
\mbox{\epsfig{file=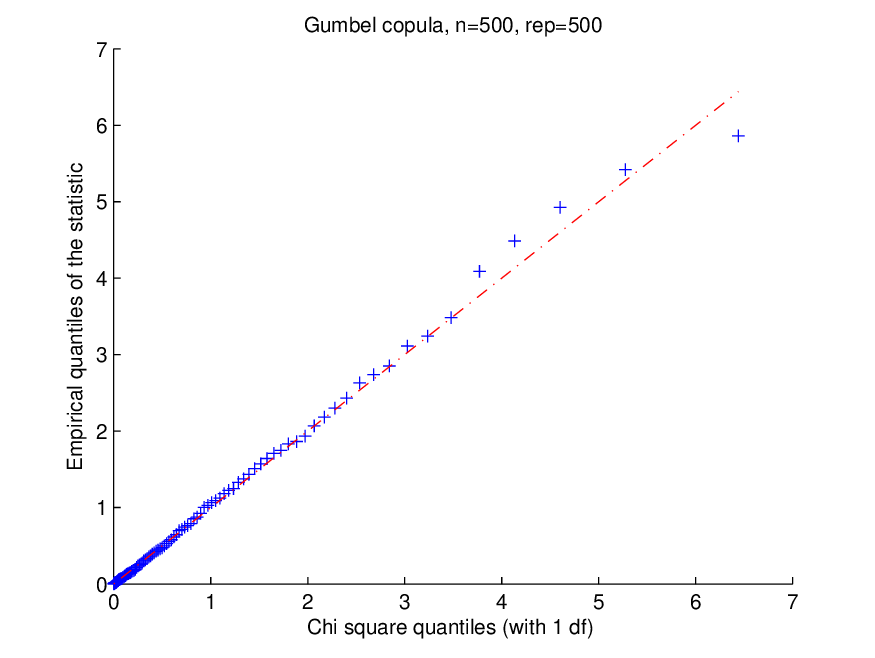,width=5cm,height=4cm}}
\end{tabular}
\end{center}
\end{boxedminipage}
\end{center}
\caption{The Q-Q Plots  of the statistic ${\bf
T}_n$ versus its limit law $\chi^2_1$}
\end{figure}
\subsection{Power comparison}
Figure 2 graphically compares the power of three different tests for bivariate
independence based on the statistics $\mathbf{T}_n$,
$\tau_n$ and $\rho_n$, where $\mathbf{T}_n$ is defined in (\ref{mmm}),
$\tau_n$ is the Kendall $\tau$
statistic,
$$\tau_n:=\frac{2}{n^2 -n}\sum_{1\leq i<j \leq n}\mbox{sign}
(R_{i}-R_{j})\mbox{sign}(S_i-S_j).$$ Here, $(R_{1},S_{1}),
\ldots,(R_{n} ,S_{n})$ are the paired rank statistics pertaining
to the sample
$(X_{1,1},X_{2,1}),\\ \ldots ,(X_{1,n},X_{2,n})$. $\rho_n$ is the Spearman $\rho$ statistic, defined by
$$\rho_n:=-3\frac{n+1}{n-1}+\frac{12}{n^3-n}\sum_{i=1}^n R_iS_i.$$
We consider the FGM copula and the marginals are exponential with
parameter one. The power functions of the three statistics are
plotted as a function of $\theta \in [0, 0.5]$ for sample sizes
$n=100$ and $n=500$. Each power entry was obtained from $1000$
independent runs. Looking at Figure 2, we can see  that the test
based on $\mathbf{T}_n$ is superior to $\tau_n$ and $\rho_n$. On the
right panel of Figure 2, for large samples (n=500), we see  that
the power of the tests based on Kendall's tau and $\mathbf{T}_n$
are almost equal.
%Note also that from this illustration, we can see that the proposed test $\mathbf{T}_n$ has a good performance.
\begin{figure}[!h]\label{comp1}
\hspace{-1cm} \centering
\mbox{\epsfig{file=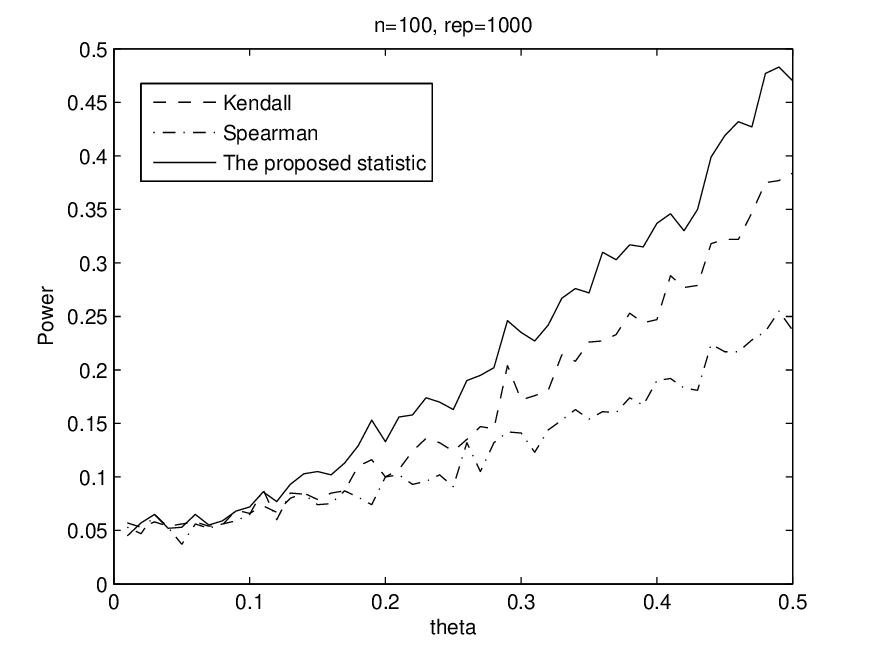,width=7cm,height=6cm}}
\mbox{\epsfig{file=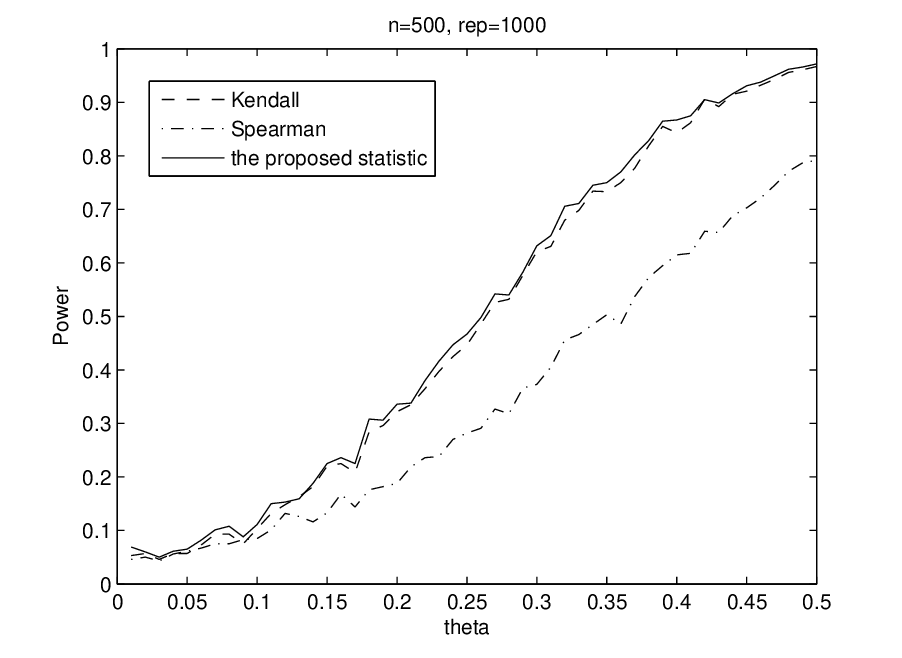,width=7cm,height=6cm}}
\caption{Power comparison}
\end{figure}
\section{Conclusion}
We have proposed  a new test procedure based on $\chi^2$-divergence and
duality technique in the framework of parametric copulas. It seems
that the procedure introduced here is particularly well adapted to
the boundary problem. For point estimation, the estimator based on $\chi^2$-divergence when we extend the parameter space, may not have a meaningful
interpretation and most probably has a larger
mean square error. However, Theorem \ref{th} may be useful to construct the confidence region for the parameter of interest in connection with the intersection method (see \cite{feng1992}) which is easy
to implement, while the asymptotic distribution
of the maximum likelihood estimator is difficult to use, since the limiting distributions are
complex when incorporating boundary constraints.

\noindent For clarity, our arguments are presented in the bivariate
context.
Note, however, that a
$d$-variate generalization is possible with obvious changes of
notation and modification assumptions. We mention that the popular multivariate Archimedean
copulas
$$
C_\theta(u_1,\ldots,u_d) = \phi^{-1} \left\{\phi(u_1) + \cdots + \phi(u_d)\right\}
$$
can only admit positive dependence, so that as long as the independence copula
belongs to a given Archimedean family, $\theta_0$ is on the boundary.
Note, also that the proposed  test of independence is not \emph{omnibus} but depends on the hypothesis that the dependence structure belong to a certain family.
\section{Appendix}\label{appendix}
\noindent First we give a technical Lemma which we will use to
prove our results.
\begin{lemma}\label{lemma1}
Let $\mathbf{F}(\cdot,\cdot)$ be continuous  and  $C(\cdot,\cdot)$ have continuous
partial derivatives. Assume that $\jmath$ is a continuous
function, with bounded variation on $I$. Then
\begin{equation}
\int_{I} \jmath(u_1,u_2)~d\left(C_{n}(u_1,u_2)-C(u_1,u_2)\right)=
O\left(n^{-1/2}(\log\log n)^{1/2}\right)~~ (a.s.).
\end{equation}
\end{lemma}\noindent{\bf{Proof of Lemma \ref{lemma1}.}} Recall that the \emph{modified
empirical copula} $C_n(\cdot,\cdot)$, is slightly different from the
\emph{empirical copula} $\mathds{C}_n(\cdot,\cdot)$, introduced by
\cite{deheuvels1979a}, and defined by
\begin{equation}\label{marre3}\mathds{C}_n(u_1,u_2)=\mathbf{F}_n\Big(F_{1n}^{-1}(u_1),
F_{2n}^{-1}(u_2)\Big)\quad\hbox{for}\quad(u_{1},u_{2})\in
(0,1)^2,\end{equation} where $F_{1n}^{-1}(\cdot)$ and
$F_{2n}^{-1}(\cdot)$ denote the empirical quantile functions,
associated with $F_{1n}(x_1)=\mathbf{F}_n(x_1,\infty)$ and
$F_{2n}(x_2)=\mathbf{F}_n(\infty,x_2)$, respectively, and defined by
\begin{equation}
F_{jn}^{-1}(t):=\inf\{x\in \mathds{R}\mid F_{jn}(x)\geq t\},\quad
j=1,2.
\end{equation}
Here,  $\mathbf{F}_n(\cdot,\cdot)$ denotes the joint empirical distribution
function, associated with the sample $\{(X_{1k},X_{2k});~
k=1,\ldots,n\},$  defined by
\begin{equation}
\mathbf{F}_{n}(x_1,x_2)=\frac{1}{n} \sum _{k=1}^{n}
\mathds{1}_{\left\{X_{1k}\leq
x_1\right\}}\mathds{1}_{\left\{X_{2k}\leq x_2\right\}},~ -\infty <
x_1,x_2< \infty.
\end{equation}
We know that $\mathds{C}_n(\cdot,\cdot)$ and $C_n(\cdot,\cdot)$ coincide on the grid
$\left\{\left(i/n,j/n\right),~1\leq i\leq j \leq n \right\}.$ The
subtle difference lies in the fact that $\mathds{C}_n(\cdot,\cdot)$ is
left-continuous with right-hand limits, whereas $C_n(\cdot,\cdot)$ on the other
hand is right continuous with left-hand limits. The difference
between $\mathds{C}_n(\cdot,\cdot)$ and $C_n(\cdot,\cdot)$, however, is small
\begin{equation}\label{app 1}
\sup_{{\bf{u}}\in
 I}\left|{\mathds{C}}_{n}({\bf{u}})-C_{n}({\bf{u}})\right|\leq
\max_{1\leq i,j\leq
n}\left|{\mathds{C}}_{n}\left(\frac{i}{n},\frac{j}{n}\right)-\mathds{C}_{n}
\left(\frac{i-1}{n},\frac{j-1}{n}\right)\right|\leq \frac{2}{n}.
\end{equation}
Using integration by parts, as in
\cite{fermanianradulovicdragan2004}, we see that
\begin{eqnarray}
\lefteqn{\sqrt{n}\int_{I}\jmath(u_1,u_2)~d(C_{n}- C)(u_1,u_2) =
\int_{I}\sqrt{n}(C_n- C)(u_1,u_2)~d\jmath(u_1,u_2)}\nonumber\\
& &-  \int_{I} \sqrt{n}(C_n- C)(u_1,1)~d\jmath(u_1,u_2)- \int_{I}
\sqrt{n}(C_n- C)(1,u_2)~d\jmath(u_1,u_2)\nonumber\\
& &-  \int_{[0,1]} \sqrt{n}(C_n(u_1,1)-u_1)~d\jmath(u_1,1)-
\int_{[0,1]} \sqrt{n}(C_n(1,u_2)-u_2)~d\jmath(1,u_2).\nonumber
\end{eqnarray}
Hence,
\begin{equation*}
\left|\sqrt{n}\int_{I}\jmath(u_1,u_2)~d(C_{n}-C)(u_1,u_2)\right|\leq
5\sqrt{n}\sup_{{\bf u}\in I}\left|(C_{n} - C)(\bf
u)\right|\int_{I}d\left|\jmath(\bf u)\right|.
\end{equation*}
From this and (\ref{app 1}), applying Theorem 3.1 in
\cite{deheuvels1979a}, we obtain the following result
$$\int_{I} \jmath(u_1,u_2)~d(C_{n}-C)(u_1,u_2)=
O\left(n^{-1/2}(\log\log n)^{1/2}\right)~~ (a.s.).$$ \hfill
$\blacksquare$\\
 \textbf{Proof of Theorem \ref{th}.}\label{thhh1}
(a) Under the assumption (C.1), a straightforward calculus yields
\begin{equation}\label{lerapportdevrais}
\int_{I}
{\textstyle\frac{\partial}{\partial\theta}}m(\theta_T,u_1,u_2)~dC_{\theta_T}(u_1,u_2)
= 0.
\end{equation}
Under the assumptions (C.1), and by applying Proposition 2.2 in
\cite{Genest_Ghoudi_Rivest1995}, we can see that, as
$n\rightarrow\infty$,
 \begin{equation}\label{normarapport}
 \int_{I} {\textstyle\frac{\partial}{\partial\theta}}m(\theta_T,u_1,u_2)
 ~dC_{n}(u_1,u_2) \longrightarrow \int_{I}
{\textstyle\frac{\partial}{\partial\theta}}m(\theta_T,u_1,u_2)~dC_{\theta_T}(u_1,u_2)
= 0,
 \end{equation}
and
\begin{equation}\label{varlim}
 \int_{I} {\textstyle\frac{\partial^2}{\partial\theta^2}}m(\theta_T,u_1,u_2)~dC_{n}(u_1,u_2)
 \longrightarrow \int_{I}
{\textstyle\frac{\partial^2}{\partial\theta^2}}m(\theta_T,u_1,u_2)~dC_{\theta_T}(u_1,u_2):=
-\Xi<0,
\end{equation}
almost surely. Now, for any $\theta = \theta_T + v n^{-1/3}$, with
$|v| \leq 1$, consider a Taylor expansion of  $\int
m(\theta,u_1,u_2)~dC_{n}(u_1,u_2)$ in $\theta$ in a neighborhood
of $\theta_T$, using (C.1) part 2, one finds
\begin{eqnarray}\label{de1}
 \lefteqn{n\int_{I} m(\theta,u_1,u_2)~dC_{n}(u_1,u_2)
-n\int_{I} m(\theta_{T},u_1,u_2)~dC_{n}(u_1,u_2)}\\
\nonumber &= & n^{2/3}v\int_{I}
{\textstyle\frac{\partial}{\partial\theta}}m(\theta_T,u_1,u_2)~dC_{n}(u_1,u_2)+n^{1/3}\frac{v^{2}}{2}\int_{I}
{\textstyle\frac{\partial^2}{\partial\theta^2}}m(\theta_T,u_1,u_2)~dC_{n}(u_1,u_2)
+O(1)~~(a.s.),
\end{eqnarray}
uniformly in $v$ with $|v|\leq 1$. On the other hand, under
condition (C.2), by Lemma \ref{lemma1}, we have
$$\int_{I} {\textstyle\frac{\partial}{\partial\theta}}m(\theta_T,u_1,u_2)dC_{n}(u_1,u_2)=
O(n^{-1/2}(\log\log n)^{1/2})\quad (a.s.).$$ Therefore, using
(\ref{varlim}) and (\ref{de1}), we obtain for any $\theta =
\theta_T + v n^{-1/3}$  with $|v|=1$,
\begin{eqnarray}
\nonumber \lefteqn{n\int_{I} m(\theta,u_1,u_2)~dC_{n}(u_1,u_2)
-n\int_{I} m(\theta_{T},u_1,u_2)~dC_{n}(u_1,u_2)}\\
& = & O(n^{1/6}(\log\log n)^{1/2})- 2^{-1}\Xi n^{1/3} +O(1)
~~(a.s.).
\end{eqnarray}
Observe that the right-hand side vanishes when $\theta=\theta_T$,
and that the left-hand side, by (\ref{varlim}), becomes negative
for all $n$ sufficiently large. Thus, by the continuity of $\theta
\mapsto \int m(\theta,u_1,u_2)~dC_{n}(u_1,u_2)$, it holds that as
$n\rightarrow\infty$, with probability one, $$\theta \mapsto \int
m(\theta,u_1,u_2)~dC_{n}(u_1,u_2)$$ reaches its maximum value at
some point $\widehat{\theta}_n$ in the interior of the interval
$B(\theta_{T},n^{-1/3})$. Therefore, the estimate $\widehat{\theta}_n$
satisfies\begin{equation} \int_{I}
{\textstyle\frac{\partial}{\partial\theta}}m(\widehat{\theta}_n,u_1,u_2)~dC_{n}(u_1,u_2)
= 0\mbox{ and   }\mid\widehat{\theta}_n -\theta_{T}\mid
=O(n^{-1/3}).\end{equation}
 \noindent (b) Making use of the first part of Theorem \ref{th},
and once more, by a Taylor expansion of
$$\int_{I} {\textstyle\frac{\partial}{\partial\theta}}
m(\widehat{\theta}_n,u_1,u_2)~dC_{n}(u_1,u_2),$$ with respect to
$\widehat{\theta}_n$, in the neighborhood of $\theta_T$, we obtain
that
\begin{eqnarray*}
\lefteqn{0 = \int_{I} {\textstyle\frac{\partial}
{\partial\theta}}m(\widehat{\theta}_n,u_1,u_2)~dC_{n}(u_1,u_2)  =
\int_{I}
{\textstyle\frac{\partial}{\partial\theta}}m(\theta_T,u_1,u_2)~dC_{n}(u_1,u_2)
}\\&& ~~~~~~~~~~~~~~+ (\widehat{\theta}_n - \theta_{T})\int_{I}
{\textstyle\frac{\partial^2}{\partial\theta^2}}m(\theta_T,u_1,u_2)~dC_{n}(u_1,u_2)
+o(n^{-1/2}).
\end{eqnarray*}
Hence,
\begin{equation}\label{equivaassmpt}
\sqrt{n}(\widehat{\theta}_n - \theta_{T}) = \left(-\int_{I}
{\textstyle\frac{\partial^2}{\partial\theta^2}}m(\theta_T,u_1,u_2)
~dC_{n}(u_1,u_2)\right)^{-1}\sqrt{n}W_n(\theta_{T}) +o(1),
\end{equation}
where $$W_n(\theta_T):=\int_{I}
{\textstyle\frac{\partial}{\partial\theta}}m(\theta_T,u_1,u_2)~
dC_n(u_1,u_2).$$ Applying Proposition 3 page 362 in
\cite{tsukahara2005}, under assumptions (C.1) part 1 and (C.3), as $n\rightarrow\infty$, we
have
\begin{equation}\label{normarlapport1}
\sqrt{n}W_{n}(\theta_{T}) \stackrel{d}{\longrightarrow}
N\left(0,\Sigma^2\right).
\end{equation}
Finally, by combining (\ref{normarlapport1}) and (\ref{varlim}) in
connection  with Slutsky's Theorem, as $n\rightarrow\infty$,we conclude that
\begin{equation}
\sqrt{n}(\widehat{\theta}_n - \theta_{T}) \stackrel{d}{\longrightarrow
}N(0,\Sigma^2/\Xi^2).
\end{equation}
This completes the proof. \hfill $\blacksquare$\vskip5pt \noindent
\textbf{Proof of Theorem \ref{thh}.}\label{thhh2} Assume that
$\theta_T=\theta_0$. From (\ref{equivaassmpt}), using
(\ref{varlim}) and (\ref{normarlapport1}), we obtain
\begin{equation}\label{eqnr 1}
\sqrt{n}\left(\widehat{\theta}_n-\theta_T\right)=-\frac{1}{\Sigma}
\sqrt{n}W_n(\theta_T)+o_P(1).
\end{equation}
Expanding in Taylor series ${\bf
T}_n(\theta_0,\widehat{\theta}_n)$ in $\widehat{\theta}_n$ around
$\theta_T$, we get
\begin{equation}\label{eqnr 2}
{\bf T}_n(\theta_0,\widehat{\theta}_n) =
2nW_n(\theta_T)(\widehat{\theta}_n-\theta_T)-\Sigma n
(\widehat{\theta}_n-\theta_T)^2+o_P(1).
\end{equation}
Now, use (\ref{eqnr 1}) combined with (\ref{eqnr 2}) to obtain
\begin{equation}\label{eqnr 3}
{\bf T}_n(\theta_0,\widehat{\theta}_n) = \frac{1}{\Sigma} n
W_n(\theta_T)^2+o_P(1).
\end{equation}
By (\ref{normarlapport1}), as $n\rightarrow\infty$, we have
\begin{equation}\label{eqnr 4}
\sqrt{n}W_{n}(\theta_{T}) \stackrel{d}{\longrightarrow} N\left(0,\Sigma^2\right)
\end{equation}
in distribution. When $\theta_T=\theta_0$, under Assumption (C.4),
we can see that $\Sigma^2$ in (\ref{eqnr 4}) is equal to $1/\Xi$;
see Proposition 2.2 in \cite{Genest_Ghoudi_Rivest1995}. Combining
this with (\ref{eqnr 3}) to conclude that \[{\bf
T}_n\stackrel{d}{\longrightarrow}\chi_{1}^2,\]  under the null
hypothesis $\mathscr{H}_0$. \hfill $\blacksquare$\vskip5pt

\noindent \noindent \textbf{Proof of Theorem \ref{thhll}.}
Rewriting $\frac{{\bf T}_n}{2n}$ as
\begin{equation}
\frac{{\bf T}_n(\theta_0,\widehat{\theta}_n)}{2n}= \sup_{\theta
\in \Theta_e} \int_{I} m(\theta,u_1,u_2)~dC_{n}(u_1,u_2),
\end{equation}
and making use of a Taylor expansion of $\frac{{\bf
T}_n(\theta_0,\widehat{\theta}_n)}{2n}$, with respect to
$\widehat{\theta}_n$, in a neighborhood of $\theta_T$, under (C.5)
to obtain
\begin{equation*}
\frac{{\bf T}_n(\theta_0,\widehat{\theta}_n)}{2n}= \int_{I}
m(\theta_T, u_1,u_2)~dC_n(u_1,u_2)+o_P(n^{-1/2}).
\end{equation*}
Hence, one finds
\begin{eqnarray*}
\sqrt{n}\left(\frac{{\bf
T}_n}{2n}-\chi^2(\theta_0,\theta_T)\right)
&=&\sqrt{n}\left(\int_{I} m(\theta_T,u_1,u_2)~dC_n-\int_{I}
m(\theta_T,u_1,u_2)~dC_{\theta_T}\right)\\&& +o_P(1).
\end{eqnarray*}
Finally, under condition (C.6), application once more of
Proposition 3 page 362 in \cite{tsukahara2005},  concludes the
proof. \hfill
$\blacksquare$ \\

%\noindent \textbf{Acknowledgement.} The authors are grateful to the two referees for their useful suggestions
%and constructive criticisms   on earlier drafts of this work.


\begin{thebibliography}{}

\bibitem[Bouzebda and Keziou(2008)]{bouzebda-keziou2008}
Bouzebda, S. and Keziou, A. (2008).
\newblock A test of independence in some copula models.
\newblock {\em Math. Methods Statist.}, {\bf 17}(2), 123--137.

\bibitem[Broniatowski and Keziou(2006)]{BK2005111}
Broniatowski, M. and Keziou, A. (2006).
\newblock Minimization of {$\phi$}-divergences on sets of signed measures.
\newblock {\em Studia Sci. Math. Hungar.}, {\bf 43}(4), 403--442.

\bibitem[Deheuvels(1979a)]{deheuvels1979a}
Deheuvels, P. (1979a).
\newblock La fonction de d\'ependance empirique et ses propri\'et\'es. {U}n
  test non param\'etrique d'ind\'ependance.
\newblock {\em Acad. Roy. Belg. Bull. Cl. Sci. (5)}, {\bf 65}(6), 274--292.

\bibitem[Deheuvels(1979b)]{deheuvels1979b}
Deheuvels, P. (1979b).
\newblock Propri\'et\'es d'existence et propri\'et\'es topologiques des
  fonctions de d\'ependance avec applications \`a la convergence des types pour
  des lois multivari\'ees.
\newblock {\em C. R. Acad. Sci. Paris S\'er. A-B}, {\bf 288}(2), A145--A148.

\bibitem[Deheuvels(1980)]{Deheuvels1981b}
Deheuvels, P. (1980).
\newblock Nonparametric test of independence.
\newblock In {\em Nonparametric asymptotic statistics (Proc. Conf., Rouen,
  1979) (French)}, volume 821 of {\em Lecture Notes in Math.}, pages 95--107.
  Springer, Berlin.

\bibitem[Deheuvels(1981)]{deheuvels1981}
Deheuvels, P. (1981).
\newblock A {K}olmogorov-{S}mirnov type test for independence and multivariate
  samples.
\newblock {\em Rev. Roumaine Math. Pures Appl.}, {\bf 26}(2), 213--226.

\bibitem[Feng and McCulloch(1992)]{feng1992}
Feng, Z. and McCulloch, C.~E. (1992).
\newblock Statistical inference using maximum likelihood estimation and the
  generalized likelihood ratio when the true parameter is on the boundary of
  the parameter space.
\newblock {\em Statist. Probab. Lett.}, {\bf 13}(4), 325--332.

\bibitem[Fermanian {\em et~al.}(2004a)]{fermanianradulovicdragan2004}
Fermanian, J.-D., Radulovi{\'c}, D., and Wegkamp, M. (2004a).
\newblock Weak convergence of empirical copula processes.
\newblock {\em Bernoulli}, {\bf 10}(5), 847--860.

\bibitem[Fermanian {\em et~al.}(2004b)]{fermanianradulovicdragan2004report}
Fermanian, J.-D., Radulovi{\'c}, D., and Wegkamp, M. (2004b).
\newblock Weak convergence of empirical copula processes.
\newblock {\em http://www.crest.fr/pageperso/fermanian/fermanian.htm}.

\bibitem[Galambos(1975)]{Galambos1975}
Galambos, J. (1975).
\newblock Order statistics of samples from multivariate distributions.
\newblock {\em J. Amer. Statist. Assoc.}, {\bf 70}(351, part 1), 674--680.

\bibitem[Genest {\em et~al.}(1995)]{Genest_Ghoudi_Rivest1995}
Genest, C., Ghoudi, K., and Rivest, L.-P. (1995).
\newblock A semiparametric estimation procedure of dependence parameters in
  multivariate families of distributions.
\newblock {\em Biometrika}, {\bf 82}(3), 543--552.

\bibitem[Gumbel(1960)]{Gumbel1960}
Gumbel, E.~J. (1960).
\newblock Bivariate exponential distributions.
\newblock {\em J. Amer. Statist. Assoc.}, {\bf 55}, 698--707.

\bibitem[H{\"u}sler and Reiss(1989)]{Husler_Reiss1989}
H{\"u}sler, J. and Reiss, R.-D. (1989).
\newblock Maxima of normal random vectors: between independence and complete
  dependence.
\newblock {\em Statist. Probab. Lett.}, {\bf 7}(4), 283--286.

\bibitem[Joe(1993)]{Joe1993}
Joe, H. (1993).
\newblock Parametric families of multivariate distributions with given margins.
\newblock {\em J. Multivariate Anal.}, {\bf 46}(2), 262--282.

\bibitem[Joe(1997)]{Joe1997}
Joe, H. (1997).
\newblock {\em Multivariate models and dependence concepts}, volume~73 of {\em
  Monographs on Statistics and Applied Probability}.
\newblock Chapman \& Hall, London.

\bibitem[Keziou and Leoni-Aubin(2005)]{KeziouLeoni2005}
Keziou, A. and Leoni-Aubin, S. (2005).
\newblock Test of homogeneity in semiparametric two-sample density ratio
  models.
\newblock {\em C. R. Math. Acad. Sci. Paris}, {\bf 340}(12), 905--910.

\bibitem[Keziou and Leoni-Aubin(2007)]{KeziouLeoni2007}
Keziou, A. and Leoni-Aubin, S. (2007).
\newblock On empirical likelihood for semiparametric two-sample density ratio
  models.
\newblock {\em Journal of Statistical Planning and Inference; In press}.

\bibitem[Kimeldorf and Sampson(1975a)]{Kimeldorf1075_1}
Kimeldorf, G. and Sampson, A. (1975a).
\newblock One-parameter families of bivariate distributions with fixed
  marginals.
\newblock {\em Comm. Statist.}, {\bf 4}, 293--301.

\bibitem[Kimeldorf and Sampson(1975b)]{Kimeldorf1075_2}
Kimeldorf, G. and Sampson, A. (1975b).
\newblock Uniform representations of bivariate distributions.
\newblock {\em Comm. Statist.}, {\bf 4}(7), 617--627.

\bibitem[Nelsen(1999)]{Nelsen1999}
Nelsen, R.~B. (1999).
\newblock {\em An introduction to copulas}, volume 139 of {\em Lecture Notes in
  Statistics}.
\newblock Springer-Verlag, New York.

\bibitem[Oakes(1994)]{Oakes1994}
Oakes, D. (1994).
\newblock Multivariate survival distributions.
\newblock {\em J. Nonparametr. Statist.}, {\bf 3}(3-4), 343--354.

\bibitem[Qin and Lawless(1994)]{qinlawless1994}
Qin, J. and Lawless, J. (1994).
\newblock Empirical likelihood and general estimating equations.
\newblock {\em Ann. Statist.}, {\bf 22}(1), 300--325.

\bibitem[R{\"u}schendorf(1976)]{Ruschendorf1976}
R{\"u}schendorf, L. (1976).
\newblock Asymptotic distributions of multivariate rank order statistics.
\newblock {\em Ann. Statist.}, {\bf 4}(5), 912--923.

\bibitem[Ruymgaart(1974)]{Ruymgaart1974}
Ruymgaart, F.~H. (1974).
\newblock Asymptotic normality of nonparametric tests for independence.
\newblock {\em Ann. Statist.}, {\bf 2}, 892--910.

\bibitem[Ruymgaart {\em et~al.}(1972)]{Ruymgaart_Shorack_Zwet1972}
Ruymgaart, F.~H., Shorack, G.~R., and van Zwet, W.~R. (1972).
\newblock Asymptotic normality of nonparametric tests for independence.
\newblock {\em Ann. Math. Statist.}, {\bf 43}, 1122--1135.

\bibitem[Self and Liang(1987)]{selflieng1987}
Self, S.~G. and Liang, K.-Y. (1987).
\newblock Asymptotic properties of maximum likelihood estimators and likelihood
  ratio tests under nonstandard conditions.
\newblock {\em J. Amer. Statist. Assoc.}, {\bf 82}(398), 605--610.

\bibitem[Shih and Louis(1995)]{Shih_Louis1995}
Shih, J.~H. and Louis, T.~A. (1995).
\newblock Inferences on the association parameter in copula models for
  bivariate survival data.
\newblock {\em Biometrics}, {\bf 51}(4), 1384--1399.

\bibitem[Sklar(1959)]{Sklar1959}
Sklar, M. (1959).
\newblock Fonctions de r\'epartition \`a {$n$} dimensions et leurs marges.
\newblock {\em Publ. Inst. Statist. Univ. Paris}, {\bf 8}, 229--231.

\bibitem[Tsukahara(2005)]{tsukahara2005}
Tsukahara, H. (2005).
\newblock Semiparametric estimation in copula models.
\newblock {\em Canad. J. Statist.}, {\bf 33}(3), 357--375.

\bibitem[Wang and Ding(2000)]{Wang_Ding2000}
Wang, W. and Ding, A.~A. (2000).
\newblock On assessing the association for bivariate current status data.
\newblock {\em Biometrika}, {\bf 87}(4), 879--893.

\end{thebibliography}
\end{document}